\documentclass[12pt,a4paper]{article}
\usepackage{amsfonts}
\usepackage{amsmath}

\newtheorem{theorem}{Theorem}[section]
\newtheorem{lemma}{Lemma}[section]
\newtheorem{prop}{Proposition}[section]
\newtheorem{cor}{Corollary}[section]

\newcommand{\R}{\mathbb{R}}

\newcommand{\nN}{n \in \mathbb{N}}

\newcommand{\C}{\mathbb{C}}
\newcommand{\E}{\mathbb{E}}

\newcommand{\cadlag}{c\`adl\`ag}

\newcommand{\g}{\textbf{g}}

\newcommand{\tr}{\mbox{tr}}

\newcommand{\Dom}{\mbox{Dom}}
\newcommand{\Ran}{\mbox{Ran}}

\newcommand{\bean}{\begin{eqnarray*}}
\newcommand{\eean}{\end{eqnarray*}}

\newcommand{\la}{\langle}
\newcommand{\ra}{\rangle}

\newcommand{\Z}{\mathbb{Z}_{+}}

\newcommand{\G}{\widehat{G}}

\newcommand{\cL}{{\cal L}}

\date{}

\begin{document}

\title{Pseudo Differential Operators and Markov Semigroups on Compact Lie Groups}

\author{ David Applebaum, \\ School of Mathematics and Statistics,\\ University of
Sheffield,\\ Hicks Building, Hounsfield Road,\\ Sheffield,
England, S3 7RH\\ ~~~~~~~\\e-mail: D.Applebaum@sheffield.ac.uk }

\maketitle

\begin{abstract} We extend the Ruzhansky-Turunen theory of pseudo differential operators on compact Lie groups into a tool that can be used to investigate group-valued Markov processes in the spirit of the work in Euclidean spaces of N.Jacob and collaborators. Feller semigroups, their generators and resolvents are exhibited as pseudo-differential operators and the symbols of the operators forming the semigroup are expressed in terms of the Fourier transform of the transition kernel.
The symbols are explicitly computed for some examples including the Feller processes associated to stochastic flows arising from solutions of stochastic differential equations on the group driven by L\'{e}vy processes. We study a family of L\'{e}vy-type linear operators on general Lie groups that are pseudo differential operators when the group is compact and find conditions for them to give rise to symmetric Dirichlet forms.

\vspace{5pt}

\noindent  {\it Key words and phrases.}~Feller semigroup, pseudo differential operator, symbol, Fourier transform, Peter Weyl theorem, Lie group, Lie algebra,
convolution semigroup, Courr\`{e}ge-Hunt operator, Sobolev space,  Dirichlet
form, Beurling-Deny representation.

\vspace{5pt}

\noindent {\it AMS 2000 subject classification.}~58J40, 47D07, 60B15, 60J25, 43A77

\end{abstract}

\section{Introduction}

The theory of pseudo differential operators has been a major theme
within modern analysis since the early 1960s and there are
important applications e.g. to  partial differential equations,
singular integral operators and index theory (see e.g. \cite{Tay}
and \cite{Hor} for details.) In the last twenty years the theory
has begun to interact with the study of stochastic processes,
principally through the work of Niels Jacob and his collaborators.
A key starting point is that the Feller semigroup and
infinitesimal generator of a convolution semigroup of probability
measures on $\R^{n}$ are naturally represented as pseudo
differential operators and the symbols of these operators may be
obtained from the classical L\'{e}vy-Khintchine formula. This case
is somewhat special in that the symbol is a function of the
Fourier variable only and this corresponds to the fact that the
characteristics of the measures have no state space dependence.
The key insight of Jacob was to consider pseudo differential
operators having more complicated symbols wherein the structure
given by the L\'{e}vy-Khintchine formula is retained but the
characteristics are now functions on $\R^{n}$.  The programme is
then to construct Feller-Markov processes whose associated
transition semigroups have generators coinciding with these pseudo
differential operators on a suitable domain. This was carried out
by using the Hille-Yosida-Ray theorem in \cite{Jac0} and by
solving the associated martingale problem in \cite{Hoh1, Hoh2}.
Applications to the construction of Dirichlet forms may be found
in \cite{HJ}. A systematic development of this theory can be found in
the three-volume work \cite{Jac2, Jac3, Jac4} (see also the
earlier concise and very readable account in \cite{Jac1}.)

The concept of pseudo differential operator may be extended to manifolds by using local co-ordinates, however when the manifold is a Lie group there is an alternative global formulation due to M.Ruzhansky and V.Turunen. Recall that if $f \in L^{1}(\R^{n}, \C)$ and its Fourier transform $\widehat{f}$ is also integrable then we have the Fourier inversion formula
$$ f(x) = \frac{1}{(2\pi)^{\frac{n}{2}}}\int_{\R^{n}}e^{iu \cdot x}\widehat{f}(u)du,$$
and if ${\cal A}$ is a pseudo differential operator with symbol $\sigma: \R^{n} \times \R^{n} \rightarrow \C$ then
\begin{equation} \label{psi}
 {\cal A}f(x) = \frac{1}{(2\pi)^{\frac{n}{2}}}\int_{\R^{n}}e^{iu \cdot x}\sigma(x,u)\widehat{f}(u)du.
 \end{equation}
 Now observe that $\R^{n}$ is its own dual group and that for each $u \in \R^{n}$, the
 character $e^{i(u, \cdot)}$ is an irreducible representation of $\R^{n}$. We can extend
 these concepts to arbitrary compact groups (Lie structure is not strictly needed here) by replacing
 the characters with general irreducible representations, defining the symbol $\sigma$  of
 a pseudo differential operator to be a suitable mapping defined on $G \times \G$ (where $\G$ is
 the set of all irreducible representations of $G$ modulo unitary equivalence) and
 utilising the Fourier inversion formula of Peter-Weyl theory (see section 3 below for
 details.) Note that for each $g \in G, \pi \in \G, \sigma(g, \pi)$ will be a matrix
 acting on the same finite-dimensional complex vector space as $\pi(g)$. A detailed
 development of this theory can be found in the monograph \cite{RT}, see also \cite{RTW}
 for recent developments including applications to global hypoellipticity.

 The purpose of the current paper is to begin the programme of developing the
  Ruzhansky-Turunen theory of pseudo differential operators on compact groups into a tool
  to study Markov processes in the spirit of Jacob's work on Euclidean space. In fact this
  task was started in the recent paper \cite{App2} but we only worked therein with
  convolution semigroups of measures (i.e. group-valued L\'{e}vy
  processes.) In the case where the group is abelian, there was
  also some work in this direction during the 1980s by N.Jacob
  \cite{Jacold}.
  Our first observation is that when we move to more general Markov processes we need to
  slightly extend the definition of pseudo differential operator from \cite{RT} where it
  is assumed that the operator preserves smooth functions. Essentially we will work in an $L^{2}$ rather than a $C^{\infty}$ framework.
  We motivate and present the relevant definitions in section 3 of the paper and also take the opportunity to review key ideas and concepts from \cite{RT}.

 The organisation of the paper is as follows. In section 2 we collect together some useful results about Feller semigroups in general locally compact spaces. As discussed above section 3 contains the definition and basic properties we'll need of general pseudo differential operators on compact groups. In section 3 we show that every Feller semigroup consists of pseudo differential operators and that these can be expressed in terms of the Fourier transform of the transition probability kernel of the associated process. In \cite{App2} (see also \cite{App1}) it was shown that the symbols of a convolution semigroup themselves form a contraction semigroup of matrices. Furthermore the symbol of the generator of the Feller semigroup is precisely the generator of the semigroup formed by the symbols.  We show here that the convolution structure is also a necessary condition for the symbols to enjoy the semigroup property. Nonetheless in the more general case, there is still sufficient regularity for us to be able to capture the symbol of the generator by differentiating that of the semigroup at the origin.

 In section 4, we collect together a number of examples of Feller processes on compact groups where we can calculate the symbol explicitly. These include Feller's pseudo-Poisson process, semigroups associated to stochastic flows obtained by solving stochastic differential equations on the group that are driven by L\'{e}vy processes and processes that are obtained by subordination. Finally in section 5, we study a general class of L\'{e}vy-type linear operators to which we can associate a symbol and thus consider as pseudo differential operators. We call these {\it Courr\`{e}ge-Hunt operators} as they are obtained by replacing the characteristics in Hunt's celebrated formula \cite{Hu} for the infinitesimal generator of a convolution semigroup with functions on the group in the spirit of Courr\`{e}ge's general form of linear operators which satisfy the positive maximal principle (\cite{Courr, Courr1}). Note however that the local part of the operator in \cite{Courr1} is expressed in local co-ordinates while we employ the Lie algebra of the group to work globally (see \cite{App4} for an earlier less general approach to this problem.) We may also regard Courr\`{e}ge-Hunt operators as natural generalisations to Lie groups of of the L\'{e}vy-type operators first studied by Komatsu \cite{Kom} and Stroock \cite{Stro} in the 1970s.

  We analyse these operators in general Lie groups (no compactness assumption is necessary at this stage) and show that (under technical conditions) they are bounded operators from the second order Sobolev space to $L^{2}(G)$ (the $L^{2}$ space of Haar measure.). Under further technical conditions, we are able to explicitly compute the adjoint and find sufficient conditions for symmetry that naturally generalise the case where the operator is the generator of a convolution semigroup (see \cite{App2}, \cite{Kun}.) Finally we are able to associate a symmetric Dirichlet form to the operator and hence realise it as the (restriction of) the generator of a Hunt process. We emphasise that the conditions for symmetry in $L^{2}(G)$ are rather strong. For future work on symmetric Courr\`{e}ge-Hunt operators it may be worthwhile to seek infinitesimal invariant measures (see \cite{ARW}) and work in the corresponding $L^{2}$ space.


 \vspace{5pt}

{\it Notation.} Let $E$ be a locally compact topological space.
The Borel $\sigma$-algebra of $E$ will be denoted by ${\cal
B}(E)$. We write $B_{b}(E)$ for the real linear space of all bounded
Borel measurable real-valued functions on $E$. It becomes a Banach
space under the supremum norm $||f||_{\infty}:=\sup_{x \in
E}|f(x)|$. The closed subspace of $B_{b}(E)$ comprising continuous functions will be denoted $C_{b}(E)$  and $C_{0}(E)$ will designate the closed subspace of continuous functions which vanish at infinity. The identity operator on $C_{0}(E)$ is denoted by $I$. In the case where $E$ is a finite dimensional real manifold, $C^{k}(E)$ will denote the linear space of $k$-times continuously differentiable functions on $E$. The smooth functions of compact support on $E$ are denoted by $C_{c}^{\infty}(E)$ and these are a dense linear manifold in $C_{0}(E)$. We will need both real and complex functions spaces in this paper. If no field is explicitly indicated you should assume that we are working in a real space. A real valued function from $[0, \infty)$ to $E$ is said to be \cadlag~if it is right continuous and left limits exist at each point.

 If $A$ is a set in a topological space then $\overline{A}$ denotes it
closure. Similarly $\overline{T}$ will denote the closure of a
closable linear operator $T$ defined on a linear submanifold of a Banach space. If $V$ is
a complex Banach space then ${\cal L}(V)$ will denote the algebra
of all continuous linear operators on $V$. We will employ the
Einstein summation convention where appropriate.

\section{Feller Processes and Semigroups}

Let $(\Omega, {\cal F}, P)$ be a probability space wherein the
$\sigma$-algebra ${\cal F}$ is equipped with a filtration $({\cal
F}_{t}, t \geq 0)$ of sub-$\sigma$-algebras. Let $Z = (Z(t), t
\geq 0)$ be a time-homogeneous Markov process with respect to the
given filtration with state space a locally compact topological
space $E$. We will denote the random variable $Z(t)$ by $Z_{x}(t)$
when we are given the initial datum $Z_{0} = x$ (a.s.) where $x \in E$. We recall
that the transition probabilities of the process are defined by
$p_{t}(x, A):= P(Z(t) \in A | Z(0) = x)$ for each $t \geq 0, x \in
E, A \in {\cal B}(E)$ and we assume that these are such that the
mappings $x \rightarrow p_{t}(x, A)$ are measurable for each $t$
and $A$. We then get a one-parameter contraction semigroup
$(T_{t}, t \geq 0)$ on $B_{b}(E)$ via the prescription
\bean  T_{t}f(x) & = & \E(f(Z(t))|Z(0) = x) \\
& = & \int_{E}f(y)p_{t}(x,dy). \eean This semigroup is also
conservative in the sense that $T_{t}1 = 1$ for each $t \geq 0$.
We say that $(T_{t}, t \geq 0)$ is a {\it Feller semigroup} (in
which case $Z$ is said to be a Feller process) if
\begin{enumerate} \item[(F1)] $T_{t}(C_{0}(E)) \subseteq C_{0}(E)$
for each $t \geq 0$, \item[(F2)] $\lim_{t \rightarrow 0}||T_{t}f -
f||_{\infty} = 0$.
\end{enumerate}
In this case $(T_{t}, t \geq 0)$ is a strongly continuous
contraction semigroup on $C_{0}(E)$ and we denote its
infinitesimal generator by ${\cal A}$. We will also require the
resolvent $R_{\lambda}({\cal A}): = (\lambda I - {\cal A})^{-1}$
which is a bounded operator on $C_{0}(E)$ for all $\lambda > 0$.

Note that to establish (F2) it is sufficient (by a standard
$\frac{\epsilon}{3}$-argument) to verify that the limit vanishes
for arbitrary $f$ in a dense linear subspace of $C_{0}(E)$ and we
will use this fact in the sequel without further comment.




\vspace{5pt}

 Let $Z$ be a Markov process. We say that there is an associated {\it solution map} if
 for each $t \geq 0$ there exists a measurable mapping
  $\Phi_{t}: E \times \Omega \rightarrow \Omega$ such that $\Phi_{t}(x, \omega) =
  Z_{x}(t)(\omega)$ for all $x \in E, \omega \in \Omega$. This structure typically arises
  when $Z$ is the solution of a stochastic differential equation (SDE) wherein the
  coefficients are sufficiently regular for the solution to generate a stochastic flow of
  homeomorphisms. The following result is surely well-known but we include a brief proof
  to make the paper more complete. We will need it in section 5 of the paper when we consider examples in compact groups.

\begin{theorem} \label{MP}
Suppose that $E$ is compact. Let $Z = (Z(t), t \geq 0)$ be a
homogeneous Markov process for which there is an associated
solution map which is such that
\begin{enumerate} \item The mapping $\Phi_{t}:E \rightarrow E$ is continuous (a.s) for all $t \geq 0$.
\item $Z$ solves the ``martingale problem'' in that there exists a
densely defined linear operator $A$ in $C_{0}(E)$ with domain $D_{A}$ such that
for all $f \in D_{A}, x \in E, t \geq 0$
\begin{equation} \label{MP1}
 f(\Phi_{t}(x)) - f(x) - \int_{0}^{t}Af(\Phi_{s}(x))ds~\mbox{is a centred martingale}.
 \end{equation}
\end{enumerate}
Then $(Z(t), t \geq 0)$ is a Feller process and $A$ extends to the infinitesimal generator of the associated Feller semigroup.
\end{theorem}

{\it Proof.} Let $(T_t, t \geq 0)$ be the semigroup associated to
the process $Z$ acting on $B_{b}(E)$. Let $(x_{n}, \nN)$ be a
sequence in $E$ which converges to $x$. Then for any $f \in
C_{0}(E)$,
$$ \lim_{n \rightarrow 0}|T_{t}f(x) - T_{t}f(x_{n})| = \lim_{n \rightarrow 0}|\E(f(\Phi_{t}(x))) - f(\Phi_{t}(x_{n})))| = 0,$$
by a straightforward application of dominated convergence and so
(F1) is satisfied. For (F2) it is sufficient to take expectations
in (\ref{MP1}) to obtain, for all $f \in D_{A}, x \in E, t \geq 0$:
$$ T_{t}f(x) - f(x) = \int_{0}^{t}T_{s}Af(x)ds,$$
and the result follows easily.  $\hfill \Box$

\vspace{5pt}

From now on we will use the term {\it Feller semigroup} to mean
any strongly continuous contraction semigroup of positive
operators on $C_{0}(E)$.

 Now let $B$ be a linear operator on $C_{0}(E)$ with domain $D_{B}$. We say that
 $B$ satisfies the {\it positive maximum principle} if $f \in D_{B}$ and $f(x_{0}) =
 \sup_{x \in E}f(x) \geq 0$ implies $Bf(x_{0}) \leq 0$. It is well-known and easily
 verified that if $A$ is the infinitesimal generator of a Feller semigroup then it
 satisfies the positive maximum principle. Conversely we have

 \begin{theorem}[Hille-Yosida-Ray] \label{hyr}
 A linear operator $A$ defined on $C_{0}(E)$ with domain $D_{A}$ is closable and its
 closure is the infinitesimal generator of a Feller semigroup if and only if
 \begin{enumerate} \item[(a)] $D_{A}$ is dense in $C_{0}(E)$.
 \item[(b)] $A$ satisfies the positive maximum principle.
 \item[(c)] $\Ran(\lambda I - A)$ is dense in $C_{0}(E)$ for some $\lambda > 0$.
 \end{enumerate}
 \end{theorem}

 See e.g. \cite{EK} pp.165-6 for a proof of this result.

 Explicit characterisations of linear operators that satisfy the positive maximum
 principle can be found in \cite{Courr} in the case $E = \R^{n}$, in \cite{Courr1} when $E$ is a compact manifold and in
 \cite{BCP} when $E$ is a manifold with compact boundary.

A Feller semigroup $(T_{t}, t \geq 0)$ is said to be {\it
conservative} if its action on $B_{b}(E)$ is such that $T_{t}1 =
1$ for all $t \geq 0$. A sufficient condition for this is that
$(1,0)$ lies in the bp-closure of the graph of $A$ (see e.g. Lemma
2.1 in \cite{Sch}). If $E$ is compact (as is the case for most of
this paper) then a necessary and sufficient condition is that $1
\in \Dom(A)$ and $A1=0$.

If $E$ is also separable then given any conservative Feller
semigroup we can construct an associated Markov process with given
initial distribution and having \cadlag~paths (see Theorem 2.7 in
\cite{EK}, pp 169-70).

\section{Pseudo Differential Operators on Compact Groups}

Let $G$ be a compact group with neutral element $e$ and let $\G$
be the set of equivalence classes (with respect to unitary
isomorphism) of irreducible representations of $G$. We will
often identify equivalence classes $[\pi]$ with a typical element $\pi$. So for each $g \in G,
\pi(g)$ is a unitary matrix acting on a finite dimensional complex
vector space $V_{\pi}$ having dimension $d_{\pi}$. We fix an
orthonormal basis $(e_{i}^{(\pi)}, 1 \leq i \leq d_{\pi})$ in $V_{\pi}$ and
define the co-ordinate functions $\pi_{ij}(g) = \la
\pi(g)e_{i}^{(\pi)}, e_{j}^{(\pi)} \ra$ for each $1 \leq i,j \leq
d_{\pi}, g \in G$. Let $L^{2}(G,\C):=L^{2}(G,m;\C)$ where $m$ is
normalised Haar measure on $(G, {\cal B}(G))$. The celebrated
Peter-Weyl theorem tells us that $\{d_{\pi}^{\frac{1}{2}}\pi_{ij},
1 \leq i,j \leq d_{\pi}\}$ is a complete orthonormal basis for
$L^{2}(G,\C)$. Moreover $L^{2}(G,\C) = \overline{\cal M}$ where ${\cal
M}$ is the linear span of $\{h_{\pi}(\psi, \phi), \psi, \phi \in
V_{\pi}, \pi \in \G\}$ where for each $g \in G, h_{\pi}(\psi,
\phi)(g):= \la \pi(g)\psi, \phi \ra$. Furthermore ${\cal M}$ is also dense in $C(G)$ (with the usual topology of uniform convergence.)
See e.g. \cite{Far} for details. In the sequel we will find it convenient to sometimes
extend the notation $h_{\pi}$ to situations where $\psi$ is
replaced by a function from $G$ to $V_{\pi}$ which is such that
the mappings $g \rightarrow \la \psi(g), \phi \ra$ are measurable
for each $\phi \in V_{\pi}$. In this case we define
$h_{\pi}(\psi(\cdot),\phi)$ to be the measurable function from $G$
to $\C$ for which $h_{\pi}(\psi(\cdot), \phi)(g) =  \la
\pi(g)\psi(g), \phi \ra$ for each $g \in G$.

The Fourier transform of $f \in L^{2}(G,\C)$ is defined by
$$ \widehat{f}(\pi) = \int_{G}\pi(g^{-1})f(g)dg,$$
where we use $dg$ as shorthand for $m(dg)$ and Fourier inversion
then yields
\begin{equation} \label{FI}
 f = \sum_{\pi \in \G}d_{\pi}\tr( \widehat{f}(\pi)\pi), \end{equation} in the
$L^{2}$-sense.  In particular for almost all $g \in G$
$$ f(g) = \sum_{\pi \in \G}d_{\pi}\tr( \widehat{f}(\pi)\pi(g)).$$
Fourier transforms of distributions on $G$ are defined by using
duality (see e.g. \cite{RT} p.545)

Now let ${\cal R}(\G) = \bigcup_{\pi \in \G}V_{\pi}$. Let
$C^{\infty}(G,\C)$ be equipped with its uniform topology and let
 $A: C^{\infty}(G,\C) \rightarrow C^{\infty}(G,\C)$ be a bounded linear map. The Ruzhansky-Turunen theory of pseudo differential operators is based on the following ideas (see \cite{RT}, p.552). By the Schwartz kernel theorem we may write $A$ as a convolution operator in that there exists a distribution on $G \times G$ such that for each $f \in C^{\infty}(G,\C), g \in G$
  $$ (Af)(g) = \int_{G}f(h)R(g, h^{-1}g)dh.$$ We define the {\it symbol} of $A$ to be the mapping $\sigma_{A}: G \times \G \rightarrow {\cal
R}(\G)$ defined by $\sigma_{A}(g, \pi) = \widehat{r_{g}}(\pi)$ where $r_{g}(\cdot) = R(g, \cdot)$
so that $\sigma_{A}(g, \pi) \in {\cal L}(V_{\pi})$ for
each $g \in G, \pi \in \G$.

 We then have (see Theorem 10.4.4 in \cite{RT}, pp.552-3)

 \begin{equation} \label{psde1}
Af(g) = \sum_{\pi \in \G}d_{\pi}\tr(\sigma_{A}(g,
\pi)\widehat{f}(\pi)\pi(g)), \end{equation} for all $f \in
C^{\infty}(G,\C), g \in G$.

It then follows (see \cite{RT}, Theorem 10.4.6, pp 553-4) that for
all $1 \leq i,j \leq d_{\pi}, \pi \in \G, g \in G$:

\begin{equation} \label{psde2}
\sigma_{A}(g, \pi)_{ij} =
\sum_{k=1}^{d_{\pi}}\overline{\pi_{ki}(g)}A\pi_{kj}(g),
\end{equation}

which we can write succinctly as

\begin{equation} \label{psde3}
\sigma_{A}(g, \pi) = \pi(g^{-1})A\pi(g),
\end{equation}
with the understanding that (\ref{psde3}) is a convenient
shorthand for the precise statement (\ref{psde2}). These ideas
extend naturally to the more general case where $A$ maps
$C^{\infty}(G, \C)$ to the space of all distributions on $G$ (see
\cite{RTW}.)

We want to extend the pseudo differential operator framework to
more general operators that are encountered in Markov process
theory. In \cite{App2} where we essentially only dealt with
L\'{e}vy processes, the symbols of all the operators that we
considered were ``functions'' of $\pi$ only and this allowed us to
effectively use (\ref{psde1}) as our definition but within a
$L^{2}$ framework. Now we want to discuss Feller processes, it
seems that a more convenient definition (and certainly one that is
in line with the theory of \cite{Jac1} and \cite{Jac2}) is to take
(\ref{psde3}) (or equivalently (\ref{psde2})) as our starting
point. To be precise we say that a linear operator $A$ defined on
$L^{2}(G,\C)$ is a {\it pseudo differential operator} if

\begin{enumerate}

\item[(PD1)] ${\cal M} \subseteq \Dom(A)$

\item[(PD2)] There exists a mapping $\sigma_{A}: G \times \G
\rightarrow {\cal R}(\G)$ such that $\sigma_{A}(g, \pi) \in {\cal
L}(V_{\pi})$ for each $g \in G, \pi \in \G$ with $\sigma_{A}(g,
\pi) = \pi(g^{-1})A\pi(g)$ for all $\pi \in \G, g \in G$.

\item[(PD3)] For each $\pi \in \G$, the mapping $g \rightarrow
\sigma_{A}(g, \pi)$ is weakly (or equivalently, strongly)
measurable.

\end{enumerate}

Of course we still call $\sigma_{A}$ the symbol of $A$. Note that a necessary and sufficient condition for (PD3) is that each matrix entry of $\sigma_{A}(g, \pi)$ is a measurable complex-valued function of $g \in G$.

Note that its symbol determines the action of a pseudo differential operator uniquely on ${\cal M}$. To see this suppose that $A_{1}$ and $A_{2}$ are densely defined linear operators in $L^{2}(G)$ with ${\cal M} \subseteq \Dom(A_{1}) \cap \Dom(A_{2})$ and that $\sigma_{A_{1}}(g, \pi) = \sigma_{A_{2}}(g, \pi)$ for all $g \in G, \pi \in \G$. Then it follows from (\ref{psde2}) that $A_{1}\pi_{ij} = A_{2}\pi_{ij}$ for all $ 1 \leq i,j \leq d_{\pi}, \pi \in \G$ and hence by linearity, $A_{1}\psi = A_{2}\psi$ for all $\psi \in {\cal M}$.

We can recover the Fourier inversion representation (\ref{psde1}) in our theory as follows.

\begin{prop} Let $A$ be a pseudo differential
operator on $L^{2}(G,\C)$ with symbol $\sigma_{A}$. Suppose that $A$ is closed and that $\sum_{\pi \in \G}d_{\pi}\tr(\widehat{f}(\pi)A\pi)$ converges in $L^{2}(G, \C)$  for all $f \in \Dom(A)$. Then
\begin{equation} \label{FIO} Af(g) = \sum_{\pi \in \G}d_{\pi}\tr(\sigma_{A}(g,
\pi)\widehat{f}(\pi)\pi(g)), \end{equation} for almost all $g \in G$.
\end{prop}

{\it Proof.} By Fourier inversion (\ref{FI}) we have for $f \in \Dom(A)$
$$ f = \sum_{\pi \in \G}d_{\pi}\sum_{i,j=1}^{d_{\pi}}\widehat{f}(\pi)_{ij}\pi_{ji}.$$  It follows from the fact that $A$ is closed that for almost all $g \in G$,
 \bean Af(g) & = & \sum_{\pi \in \G}d_{\pi}\sum_{i,j=1}^{d_{\pi}}\widehat{f}(\pi)_{ij}A\pi_{ji}(g)\\
 & = & \sum_{\pi \in \G}d_{\pi}\tr( \widehat{f}(\pi)A\pi(g))\\
& = & \sum_{\pi \in \G}d_{\pi}\tr( \widehat{f}(\pi)\pi(g)\sigma_{A}(g,\pi))\\
& = & \sum_{\pi \in \G}d_{\pi}\tr(\sigma_{A}(g,
\pi)\widehat{f}(\pi)\pi(g)), \eean
where we have used (\ref{psde3}) and the fact that $\tr(XY) = \tr(YX). \hfill \Box$

\vspace{5pt}


In the sequel, when we consider Feller semigroups, we will want to
consider operators that preserve the space of continuous functions
on $G$. Indeed we will need the fact that any densely defined
linear operator on $C(G)$ (equipped with the uniform topology) is
also densely defined on $L^{2}(G)$.

\begin{prop} If $A: C(G) \rightarrow C(G)$ is a pseudo differential operator then for each $\pi \in G$ the mapping $g \rightarrow \sigma_{A}(g, \pi)$ is
weakly (equivalently, strongly) continuous.
\end{prop}

{\it Proof.} This follows immediately from (\ref{psde2}). $\hfill \Box$

\vspace{5pt}

The following result gives an equivalent characterisation of a pseudo differential
operator:

\begin{prop} \label{pseq}
A linear operator $A$ defined on $L^{2}(G,\C)$ with ${\cal M}
\subseteq \Dom(A)$ is a {\it pseudo differential operator} if and only if
there exists a mapping $\sigma_{A}: G \times \G \rightarrow {\cal
R}(\G)$ satisfying (PD3) such that $\sigma_{A}(g, \pi) \in {\cal L}(V_{\pi})$ and
\begin{equation} \label{psde4}
A h_{\pi}(\psi, \phi)(g) =  h_{\pi}(\sigma_{A}(\cdot, \pi) \psi,
\phi)(g),
\end{equation}
for all $h_{\pi} \in {\cal M}, g \in G$.
\end{prop}

{\it Proof.}~ With $h_{\pi},g$ as in the statement of the
proposition, by (\ref{psde3}): \bean A h_{\pi}(\psi, \phi)(g) & =
& \la A \pi(g) \psi, \phi \ra \\
& = &  \la \pi(g^{-1})A \pi(g) \psi, \pi(g^{-1}) \phi \ra \\
& = & \la \sigma_{A}(g, \pi)\psi, \pi(g^{-1}) \phi \ra \\
& = & \la \pi(g)\sigma_{A}(g, \pi)\psi,  \phi \ra \\
& = & h_{\pi}(\sigma_{A}(\cdot, \pi) \psi, \phi)(g), \eean where
we have used the fact that $\pi(g)^{*} = \pi(g^{-1})$. The
converse follows when we take $\psi$ and $\phi$ to be
$e_{i}^{(\pi)}$ and $e_{j}^{(\pi)}$ (respectively). We then find by (\ref{psde4}) that
$$ A\la \pi(\cdot)e_{i}^{(\pi)}, e_{j}^{(\pi)}\ra (g) = \la \pi(g)\sigma_{A}(g, \pi)e_{i}^{(\pi)}, e_{j}^{(\pi)} \ra$$
and so
$$ A\pi_{ij}(g) = \sum_{k=1}^{d_{\pi}}\pi_{ik}(g)\sigma_{A}(g, \pi)_{kj}.$$
By unitarity of the matrix $\pi(g)$ for each $1 \leq l \leq d_{\pi}$, we obtain
$$ \sum_{i=1}^{d_{\pi}}\overline{\pi_{il}(g)}A\pi_{ij}(g) = \sigma_{A}(g, \pi)_{lj}, $$
and since this holds for each $1 \leq j \leq d_{\pi}$ then (\ref{psde2}) is obtained, as required.  $~~~~~~~~~~~~~~~~~~~~~~~~~~~~~~~~~~~~~~~~~~~~~~~~~~~~~~~~~~~~~~~~~~~~~~~~~~~~~~~~~~~~~~~~~~~~~~~~~~~~\hfill \Box$

\vspace{5pt}

We remark that Proposition \ref{pseq} enables us to extend the
definition of pseudo differential operators to more general
topological groups. We will not pursue that theme further in this
article.

\section{Pseudo Differential operators and Feller semigroups}

 We begin this section with some preliminaries about measures on groups and their Fourier transforms. Recall (see e.g. \cite{Hey}, \cite{Sieb}) that if $\mu$ is  a  Borel probability measure defined on $G$ then its Fourier transform $(\widehat{\mu}(\pi), \pi \in \G)$ is defined
 by the matrix-valued integral\footnote{Note that we employ the ``probabilist's convention'' for Fourier transforms of measures and the ``analyst's convention'' for those of functions. Readers should be reassured that
this does not result in any mathematical inconsistency.}:
$$ \widehat{\mu}(\pi) = \int_{G}\pi(g)\mu(dg),$$ so that each $\widehat{\mu}(\pi)$ acts as a contraction on $V_{\pi}$.
Furthermore $\widehat{\mu}$ determines $\mu$ uniquely. If $\mu$ and $\nu$ are Borel probability measures on $G$ then for all $\pi \in \G$
$$ \widehat{\mu * \nu}(\pi) = \widehat{\mu}(\pi)\widehat{\nu}(\pi),$$
where $\mu * \nu$ denotes the convolution of $\mu$ and $\nu$, i.e. the unique Borel probability measure on $G$ for which
$$ \int_{G}f(\rho)(\mu * \nu)(d\rho) = \int_{G}\int_{G}f(\rho\tau)\mu(d\rho)\nu(d\tau),$$ for all $f \in B_{b}(G)$.

Let $(T_{t}, t \geq 0)$ be a Feller semigroup defined on $C(G)$
where $G$ is a compact group. Then for each $t \geq 0, T_{t}$ is a
densely defined (not necessarily bounded) linear operator in
$L^{2}(G)$ with ${\cal M} \subseteq \Dom(T_{t})$. Recall that for each $t \geq 0, g \in G$ the transition probability $p_{t}(g, \cdot)$ is a Borel probability measure on $G$. We denote its Fourier transform at $\pi \in \G$  by $\widehat{p_{t}}(g, \pi)$.

\begin{lemma} \label{FPD1}
For each $t \geq 0, T_{t}$ is a pseudo differential operator with symbol
\begin{equation} \label{FPD2}
                  \sigma_{t}(g, \pi) = \pi(g^{-1})\widehat{p_{t}}(g, \pi),
\end{equation}
for $g \in G, \pi \in \G$.
\end{lemma}

{\it Proof.}  From the discussion above it follows that (PD1) is satisfied. To establish the identity (\ref{FPD2}) we use
\bean \sigma_{t}(g, \pi) & = & \pi(g^{-1})T_{t}\pi(g)\\
& = & \pi(g^{-1})\int_{G}\pi(\tau)p_{t}(g , d\tau)\\
& = & \pi(g^{-1})\widehat{p_{t}}(g, \pi), \eean as required.
(PD3) follows from the facts that $g \rightarrow \widehat{p_{t}}(g, \pi)$ is weakly measurable and $g \rightarrow \pi(g)$ is strongly continuous. $\hfill \Box$

\vspace{5pt}

Note that we always have $\sigma_{0}(g, \pi) = I_{\pi}$, where $I_{\pi}$ is the identity matrix on $V_{\pi}$.

\vspace{5pt}

Example. {\it Convolution Semigroups}. Let $(\mu_{t}, t \geq 0)$ to be a vaguely continuous convolution semigroup
 of probability measures on $G$. Such semigroups of measures arise as the laws of $G$-valued L\'{e}vy processes (see e.g. \cite{Liao}). In this case for each $t \geq 0, g \in G, f \in C(G), T_{t}f(g) = \int_{G}f(g \tau)\mu_{t}(d\tau)$ and $p_{t}(g, A) = \mu_{t}(g^{-1}A)$ for each $g \in G, A \in {\cal B}(G)$. So by Lemma \ref{FPD1}
 we obtain  \bean \sigma_{t}(g, \pi) & = & \pi(g^{-1})\int_{G}\pi(\tau)\mu_{t}(g^{-1}d\tau)\\
 & = & \widehat{\mu_{t}}(\pi) \eean for all $g \in G, \pi \in \G$.
 It is shown in \cite{App1} that for each $\pi \in G, (\widehat{\mu_{t}}(\pi), t \geq 0)$ is a strongly continuous contraction semigroup of matrices acting on $V_{\pi}$. The fact that the symbol is independent of $g \in G$ is clearly related to the translation invariance of the semigroup, i.e. the fact that $L_{\rho}T_{t} = T_{t}L_{\rho}$ for all $\rho \in G, t \geq 0$ where $L_{\rho}f(g) = f(\rho^{-1}g)$ for each $f \in C(G), g \in G$ (c.f. \cite{Hu}.)

 \vspace{5pt}

 We will investigate more examples in the next section. First we establish some general properties.

 \begin{theorem} \label{gen}

 Let $(\sigma_{t}, t \geq 0)$ be the symbol of a Feller semigroup $(T_{t}, t \geq 0)$ with infinitesimal generator ${\cal A}$ and assume that ${\cal M} \subseteq \mbox{Dom}({\cal A})$. Then for each $g \in G, \pi \in \G$,

 \begin{enumerate}
 \item The mapping $t \rightarrow \sigma_{t}(g, \pi)$ is strongly differentiable in $V_{\pi}$.

 \item For each $s, t \geq 0$,
 \begin{equation} \label{notsemi}
 \sigma_{s+t}(g, \pi) = \pi(g^{-1})T_{s}\widehat{p_{t}}(\cdot, \pi)(g).
 \end{equation}

 \end{enumerate}
 \end{theorem}

 {\it Proof}. \begin{enumerate} \item  Since $(T_{t}, t \geq 0)$ is strongly differentiable on Dom$({\cal A})$ and $\pi_{ij} \in \mbox{Dom}({\cal A})$ for each $1 \leq i,j \leq d_{\pi}$, it follows by (\ref{psde2}) that each matrix entry $\sigma_{t}(g, \pi)_{ij}$ is a differentiable function of $t$. We define a $d_{\pi} \times d_{\pi}$ matrix $j_{t}(g, \pi)$ by the prescription
 $$ j_{t}(g, \pi)_{ij} = \frac{d}{dt}\sigma_{t}(g, \pi)_{ij} = \sum_{k=1}^{d_{\pi}}\overline{\pi_{ki}(g)}T_{t}{\cal A}\pi_{kj}(g).$$
 Now let $\psi = (\psi_{1}, \ldots, \psi_{d_{\pi}}) \in V_{\pi}$. Using the Cauchy-Schwarz inequality (twice) and the fact that $|\pi_{ij}(g)| \leq 1$ for each $1 \leq i,j \leq d_{\pi}, g \in G$ we have for each $t \geq 0$
  \bean & & \limsup_{h \rightarrow 0} \frac{1}{h}\left|\left|\sigma_{t+h}(g, \pi)\psi - \sigma_{t}(g, \pi)\psi - hj_{t}(g, \pi)\psi\right|\right|^{2}_{V_{\pi}}\\
 & = & \limsup_{h \rightarrow 0} \frac{1}{h}\sum_{i=1}^{d_{\pi}}\left|\sum_{j=1}^{d_{\pi}}\sum_{k=1}^{d_{\pi}}\overline{\pi_{ki}(g)}(T_{t+h}\pi_{kj}(g) - T_{t}\pi_{kj}(g) - hT_{t}{\cal A}\pi_{kj}(g))\psi_{j}\right|^{2}\\
 & \leq & d_{\pi}||\psi||^{2}_{V_{\pi}}\limsup_{h \rightarrow 0} \frac{1}{h}\sum_{i,j,k=1}^{d_{\pi}}|T_{t+h}\pi_{kj}(g) - T_{t}\pi_{kj}(g) - hT_{t}{\cal A}\pi_{kj}(g)|^{2} = 0,\eean
 and the result follows since $\pi_{kj} \in \Dom({\cal A})$.

 \item Using the Chapman-Kolmogorov equations we find that
 \bean \sigma_{s+t}(g, \pi) & = & \pi(g^{-1})T_{s+t}\pi(g)\\
 & = & \int_{G}\pi(g^{-1}\tau)p_{s+t}(g, d\tau)\\
 & = & \int_{G}\pi(g^{-1}\rho)\left(\int_{G}\pi(\rho^{-1}\tau)p_{t}(\rho, d\tau)\right)p_{s}(g, d\rho)\\
 & = & \int_{G}\pi(g^{-1}\rho)\sigma_{t}(\rho, \pi)p_{s}(g, d \rho)\\
 & = & \pi(g^{-1})T_{s}(\pi(\cdot)\sigma_{t}(\cdot, \pi))(g)\\
 & = & \pi(g^{-1})T_{s}\widehat{p_{t}}(\cdot, \pi)(g), \eean
 where we used (\ref{FPD2}) to obtain the last line.

 \end{enumerate}
 $\hfill \Box$

 For each $g \in G, \pi \in \G$ we define a $d_{\pi} \times d_{\pi}$ matrix $j(g, \pi)$ by the prescription
 \begin{equation} \label{gend}
 j(g, \pi)_{ij}: = \left.\frac{d}{dt}\sigma_{t}(g, \pi)_{ij}\right|_{t=0}
 \end{equation}
for each $1 \leq i,j \leq d_{\pi}$ (of course $j$ is just $j_{0}$ from the proof of Theorem \ref{gen}.)

 \begin{theorem} \label{JA}
 \begin{enumerate}
 \item If ${\cal M} \subseteq \Dom({\cal A})$ then ${\cal A}$ is a pseudo differential operator with symbol $j(g, \pi)$ at $\pi \in \G, g \in G$.
 \item For each $\lambda > 0, R_{\lambda}(A)$ is a pseudo differential operator with symbol $\int_{0}^{\infty}e^{-\lambda t}\sigma_{t}(g, \pi)dt$ at $\pi \in \G, g \in G$.
 \end{enumerate}
 \end{theorem}

 {\it Proof.} \begin{enumerate}
 \item Using (\ref{psde3}) we compute
 \bean \pi(g^{-1}){\cal A}\pi(g) & = & \pi(g^{-1})\left.\frac{d}{dt}T_{t}\pi(g)\right|_{t=0}\\
 & = & \left.\frac{d}{dt}\pi(g^{-1})T_{t}\pi(g)\right|_{t=0}\\
 & = & \left.\frac{d}{dt}\sigma_{t}(g, \pi)\right|_{t=0} = j(g, \pi). \eean

\item This follows from (\ref{psde3}) using the fact that $R_{\lambda}f = \int_{0}^{\infty}e^{-\lambda t}T_{t}fdt$ for each $f \in C(G). \hfill \Box$

\end{enumerate}

\vspace{5pt}

{\bf Note.} If $(T_{t}, t \geq 0)$ is the semigroup associated to
the $G$-valued Feller process $(X(t), t \geq 0)$ then Theorem
\ref{JA}(1) tells us that for all $g \in G, \pi \in \G$,
$$ j(g, \pi) = \lim_{t \rightarrow
0}\frac{1}{t}(\pi(g^{-1})\E(\pi(X(t))|X(0)=g) - I_{\pi}),$$ and
this should be compared with the corresponding analysis in
Euclidean space - see e.g. equation (0.3) on p.xx in \cite{Jac2}
and Definition 3.9.1 on p.148 of \cite{Jac4}.

\vspace{5pt}

We close this section by returning to the question of when the matrices $(\sigma_{t}(e, \pi), t \geq 0)$ form a semigroup on $V_{\pi}$. We know from (\ref{notsemi}) that this is unlikely to hold in general, but that this is the case when the transition probabilities are translates of a vaguely continuous convolution semigroup of probability measures. The next result shows that this is the only possibility.

\begin{prop} For each $t \geq 0$, let $\sigma_{t}(g, \pi)$ be the symbol of a Feller semigroup at $g \in G, \pi \in \G$. We have $\sigma_{s+t}(e, \pi) = \sigma_{s}(e, \pi)\sigma_{t}(e, \pi)$ if and only if $(p_{t}(e, \cdot), t \geq 0)$ is a vaguely continuous convolution semigroup of probability measures on $G$.
\end{prop}

{\it Proof.} By (\ref{FPD2}) $\sigma_{t}(e, \pi) = \widehat{p_{t}}(e, \pi)$ and so we have $\sigma_{s+t}(e, \pi) = \sigma_{s}(e, \pi)\sigma_{t}(e, \pi)$ if and only if
$$ \widehat{p_{s+t}}(e, \pi) = \widehat{p_{s}}(e, \pi)\widehat{p_{t}}(e, \pi) = \widehat{p_{s}* p_{t}}(e, \pi)$$ for all $\pi \in \G$,
and the result follows by uniqueness of the Fourier transform. $~~~~~~~~~~~~~~~~~~~~~~~~~~~~~~~~~~~~~~~~~~~~~~~~~~~~~~~~~~~~~~~~~~~~~~~~\hfill \Box$






\section{Examples}

\begin{enumerate}

\item {\it Feller's Pseudo-Poisson Process}

Let $(S(n), n \in \Z)$ be a Markov chain taking values in $G$ with
transition kernel $q:G \times {\cal B}(G) \rightarrow [0,1]$, so
that $q(g,B) = P(S(1) \in B | S(0) = g)$ for all $g \in G, B \in
{\cal B}(G)$. Let $(N(t), t \geq 0)$ be a Poisson process of
intensity $\lambda > 0$ that is independent of $(S(n), n \in \Z)$.
Then Feller's pseudo-Poisson process $(Z(t), t \geq 0)$ is defined
by $Z(t) = S(N(t))$ for each $t \geq 0$. It is well known (see
e.g. \cite{Abook}, \cite{EK}, \cite{Jac4}) that it is a Feller
process with bounded infinitesimal generator defined for each $f
\in L^{2}(G,\C), g \in G$ by
\begin{equation} \label{fe}
{\cal A}f(g) = \lambda \int_{G}(f(\tau) - f(g))q(g,d\tau).
\end{equation} Hence by a straightforward application of the definition
(PD2), we find that the symbol is given by
\begin{equation} \label{Feller}
j(g, \pi) = \lambda \int_{G}(\pi(g^{-1}\tau) - I_{\pi})q(g,d\tau)
\end{equation}
for all $\pi \in \G$.

\vspace{5pt}

For the remaining examples in this section, we will require that
$G$ be a compact $n$-dimensional Lie group with Lie algebra $\g$. Let $\exp$ be
the exponential map from $\g$ to $G$. Then for each $X \in \g, \pi
\in \G$, the skew-symmetric matrix $d\pi(X)$ is defined by
$$ d\pi(X) = \left.\frac{d}{du}\pi(\exp(uX))\right|_{u=0}.$$
Now consider the left invariant vector field $g \rightarrow X(g)$ for which $X(e) = X$.
It is a pseudo differential operator and its symbol is given by
\bean \sigma_{X}(g, \pi) & = & \pi(g^{-1})X\pi(g)\\
& = & \pi(g^{-1})\left.\frac{d}{du}\pi(g\exp(uX))\right|_{u=0}\\
& = & \left.\frac{d}{du}\pi(\exp(uX))\right|_{u=0}\\
& = & d\pi(X), \eean
for each $g \in G, \pi \in \G$. If $X^{\prime}$ is the right invariant vector field for which $X^{\prime}(e) = X$ then a similar argument yields
 $\sigma_{X^{\prime}}(g, \pi) = \pi(g^{-1})d\pi(X)\pi(g)$. We will make use of these facts in the remaining examples.

\item{\it Convolution Semigroups}

Let $(\mu_{t}, t \geq 0)$ be a vaguely continuous convolution
semigroup of measures on $G$ with $\mu_{0} = \delta_{e}$. As
discussed above the associated Feller semigroup on $C(G)$ is given
by $T_{t}f(g) = \int_{G}f(g\tau)\mu_{t}(d\tau)$ for $f \in C(G), t
\geq 0, g \in G$. The infinitesimal generator ${\cal A}$ was first
studied by Hunt \cite{Hu} - see also \cite{He1} and \cite{Liao}.
Let $\{X_{1}, \ldots, X_{n}\}$ a fixed basis of $\g$ and define
$$C_{2}(G):=\{f \in C(G), X_{i}f \in C(G), X_{j}X_{k} \in C(G), 1 \leq i,j,k \leq n\}.$$  It is easily verified that ${\cal M} \subseteq C_{2}(G)$. It is shown in \cite{Hu} that $C_{2}(G) \subseteq \Dom({\cal A})$ and for $f \in C_{2}(G)$ and
\begin{eqnarray} \label{Hu}
{\cal A}f(g) & = & b^{i}X_{i}f(g) + a^{ij}X_{i}X_{j}f(g) \nonumber \\
& + & \int_{G-\{e\}}[f(g\tau) - f(g) - x^{i}(\tau)X_{i}f(g)] \nu(d\tau)
\end{eqnarray}

where $b \in \R^{n}, (a_{ij})$ is a non-negative definite symmetric $n \times n$ matrix
and $\nu$ is a L\'{e}vy  measure on $G-\{e\}$ i.e. a Borel measure for which

$\int_{G-\{e\}}\left(\sum_{i=1}^{n}x_{i}(\tau)^{2} \wedge 1
\right)\nu(d\tau) < \infty$. Here $x_{i} \in C^{\infty}(G)$ for $1
\leq i \leq n$ is such that $(x_{1}, \ldots, x_{n})$ are canonical
co-ordinates in a neighborhood of $e$ with $x_{i}(e) = 0$ and
$X_{i}x_{j}(e) = \delta_{ij} (1 \leq i,j \leq n.)$

We have already seen that the symbol of $T_{t}$ is $\widehat{\mu}_{t}$ for each $t \geq 0$. It follows from (\ref{gend}) that that of ${\cal A}$ is the infinitesimal generator of the matrix semigroup $(\widehat{\mu}_{t}, t \geq 0)$. This is obtained in \cite{App1} (see also \cite{Hey1}, \cite{LoNg} and the discussion in \cite{App2}) and we have for all $g \in G, \pi \in \G$
\begin{eqnarray} \label{Hun1}
j(g, \pi) & = & b^{i}d\pi(X_{i}) + a^{ij}d\pi(X_{i})d\pi(X_{j}) \nonumber \\
& + & \int_{G-\{e\}}[\pi(\tau) - I_{\pi} - x^{i}(\tau)d\pi(X_{i})] \nu(d\tau)
\end{eqnarray}

\item {\it L\'{e}vy Flows and Diffusions}

Let $Y_{1}, \ldots, Y_{p}$ be $C^{\infty}$-vector fields on $G$ and let $X_{1}, \ldots, X_{n}$ be a basis for $\g$. Since $G$ is parallelisable, we can assert the existence of $C^{\infty}$-mappings $\gamma^{j}_{i}:G \rightarrow \R$ such that for each $g \in G, Y_{i}(g) = \gamma^{j}_{i}(g)X_{j}(g)$ for each $ 1 \leq i \leq p$. Now let $(L(t), t \geq 0)$ be an $\R^{p}$-valued \cadlag~L\'{e}vy process with characteristics $(b, a, \nu)$ and consider the stochastic differential equation

\begin{equation} \label{SDE}
  d\phi(t) = Y_{i}(\phi(t-)) \diamond dL^{i}(t),
\end{equation}

where $\diamond$ denotes the Markus canonical form (see e.g. section 6.10 in \cite{Abook}.) It is shown in \cite{ATa2} pp. 233-4 that the solution map is a stochastic flow of diffeomorphisms of $G$ and it follows by Theorem \ref{MP} that $(\phi(t), t \geq 0)$ is a Feller process. The infinitesimal generator is given by
\begin{eqnarray} \label{Mark1}
{\cal A}f(g) & = & b^{i}Y_{i}f(g) + \frac{1}{2}a^{ij}Y_{i}Y_{j}f(g) \nonumber \\
& + & \int_{\R^{p}-\{0\}}(f(\xi(y)(g)) - f(g) - y^{i}Y_{i}f(g){\bf
1}_{|y| < 1})\nu(dy)\nonumber\\
& &
\end{eqnarray}
for $f \in C^{2}(G), g \in G$. Here $(\xi(y), y \in G)$ is the flow
of diffeomorphisms of $G$ obtained by solving the ordinary
differential equation
   $$ \frac{d\xi(vy)}{dv} = y^{i}Y_{i}\xi(vy).$$
In this case the symbol of ${\cal A}$ is given by
\begin{eqnarray} \label{Mark2}
j(g, \pi) & =  & b^{i}\gamma^{j}_{i}(g)d\pi(X_{j}) \nonumber \\ & + &  \frac{1}{2} a^{ij}\gamma^{k}_{i}(g)\left[\gamma^{r}_{j}(g)d\pi(X_{k})d\pi(X_{r}) +  X_{k}(g)\gamma^{r}_{j}(g)d\pi(X_{r})\right] \nonumber \\
& + & \int_{\R^{p}-\{0\}}(\pi(g^{-1}\xi(y)(g)) - I_{\pi} - y^{i}\gamma^{k}_{i}(g)d\pi(X_{k}){\bf 1}_{|y| < 1})\nu(dy), \nonumber \\ & & \end{eqnarray}
for each $\pi \in \G, g \in G$.

Note that when $\nu \equiv 0$, then $\diamond$ is a Stratonovitch differential equation and $(\phi(t), t \geq 0)$ gives rise to a diffusion process on $G$ (see \cite{BHT} for some interesting recent work on such processes.) We remark that a general form for symbols of infinitesimal generators of Feller processes that arise as solutions of SDEs driven by L\'{e}vy processes in $\R^{n}$ may be found in \cite{ScSc}.

\item {\it Pseudo-Feller Jump Diffusions}

We consider the linear operator ${\cal A}$ defined on $C^{2}(G)$
by
\begin{equation} \label{jd}
{\cal A}f(g)  =  b^{i}Y_{i}f(g) + \frac{1}{2}a^{ij}Y_{i}Y_{j}f(g)
+ \lambda \int_{G}(f(\tau) - f(g))q(g,d\tau),
\end{equation}

for $f \in C^{2}(G), g \in G$ where the $b_{i}$s, $a_{ij}$ and
$Y_{i}$s are as in Example 3 and $\lambda$ and $q$ are as in
Example 1. In fact we can  write ${\cal A} = {\cal A}_{1} + {\cal
A}_{2}$ where ${\cal A}_{1}$ is as in (\ref{Mark1}) but with $\nu
\equiv 0$ and ${\cal A}_{2}$ is as in (\ref{fe}). Since both
${\cal A}_{1}$ and ${\cal A}_{2}$ generate one parameter
contraction semigroups in $C(G)$ so does ${\cal A}$ (see Theorem
2.7 in \cite{Kato}, Chapter IX, p.501). Since both ${\cal A}_{1}$
and ${\cal A}_{2}$ generate Feller semigroups, they satisfy the
positive maximum principle and hence ${\cal A}$ also satisfies
this principle. It follows by the Hille-Yosida-Ray theorem that
the closure of ${\cal A}$ generates a Feller semigroup. We can
easily compute the symbol of the generator to be:
\begin{eqnarray} \label{jd1}
j(g, \pi) & = & b^{i}\gamma^{j}_{i}(g)d\pi(X_{j})+ \frac{1}{2} a^{ij}\gamma^{k}_{i}(g)\left[\gamma^{r}_{j}(g)d\pi(X_{k})d\pi(X_{r})\right. \nonumber \\ & + & \left. X_{k}(g)\gamma^{r}_{j}(g)d\pi(X_{r})\right]
 +  \lambda \int_{G}(\pi(g^{-1}\tau) - I_{\pi})q(g,d\tau). \nonumber \\
 & &
\end{eqnarray} for each $\pi \in \G, g \in G$.

\item {\it Subordination}

Let $(\rho_{t}^{h}, t \geq 0)$ be the law of a
subordinator with associated Bernstein function $h:(0, \infty)
\rightarrow [0, \infty)$ so that $(\rho_{t}^{h}, t \geq 0)$ is a
vaguely continuous convolution semigroup of probability measures
on $[0, \infty)$ and for each $t \geq 0, u > 0$,
\begin{equation} \label{sub}
\int_{0}^{\infty}e^{-us}\rho_{t}^{h}(ds) = e^{-th(u)}.
\end{equation}
Note that $h$ has the generic form
$$ h(u) = bu + \int_{(0,\infty)}(1 - e^{-uy})\lambda(dy),$$
where $b \geq 0$ and $\int_{(0,\infty)}(y \wedge 1)\lambda(dy) <
\infty$ (see e.g. \cite{Sa}, section 30 and \cite{Abook}, section
1.3.2 for details.) Let $(T_{t}, t \geq 0)$ be a Feller semigroup
on $C(G)$ for which ${\cal M} \subseteq \mbox{Dom}({\cal A})$. By
Lemma \ref{FPD1} and Theorem \ref{JA}, for each $t \geq 0$ $T_{t}$
and ${\cal A}$ are pseudo differential operators and we denote
their symbols by $\sigma_{t}$ and $j$ (respectively). Appealing to
Phillip's theorem (see \cite{Phi}, \cite{Sa} pp. 212-7) we see
that subordination yields another Feller semigroup $(T_{t}^{h}, t
\geq 0)$ on $C(G)$ with infinitesimal generator ${\cal A}^{h}$
having Dom$({\cal A})$ as a core, where for each $f \in C(G)$ we
have (in the sense of Bochner integrals)
$$ T_{t}^{h}f = \int_{0}^{\infty} T_{s}f \rho_{t}^{h}(ds),$$
and for each $f \in \mbox{Dom}({\cal A})$,
$$ {\cal A}^{h}f = b{\cal A}f + \int_{(0, \infty)} (T_{s}f - f)\lambda^{h}(ds).$$
It then follows easily that ${\cal A}^{h}$ is a pseudo differential operator and its symbol is given for $g \in G, \pi \in \G$ by
$$ j^{h}(g, \pi) = bj(g, \pi) + \int_{(0, \infty)}(\sigma_{s}(g, \pi) - I_{\pi})\lambda^{h}(ds).$$

\item {\it Class $O$ Semigroups}

We have already remarked that the Feller semigroups associated with convolution semigroups of measures commute with left translation on the group. In \cite{App4} a class of Feller semigroups was investigated under which there was an obstruction to these commutation relations. To be precise we require that there exists a measurable mapping $h:\R^{+} \times G \times G \times G \rightarrow \R$ for which
  $$ L_{\rho}T_{t}f(g) = T_{t}L_{\rho}fh(t, g, \rho, \cdot)(g) ,$$
  for all $t \geq 0, g, \rho \in G.$ Such semigroups were said to be of Class $O$
  (where the $O$ stands for ``obstruction''.) Define $C_{2}^{\prime}(G)$ in exactly the same way as $C_{2}(G)$ but with the basis of left-invariant vector fields replaced by right-invariant ones.
 Under some technical conditions on $h$ we have $C_{2}^{\prime}(G) \subseteq \Dom({\cal A})$ and (using the
same notation where appropriate as was the case in Example 2) for
$f \in C_{2}^{\prime}(G)$
\begin{eqnarray} \label{Hum}
{\cal A}f(g) & = & b^{i}(g)X_{i}^{\prime}f(g) + a^{ij}(g)X_{i}^{\prime}X_{j}^{\prime}f(g) \nonumber \\
& + & \int_{G-\{e\}}[f(\tau) - f(g) - x^{i}(\tau g^{-1})X_{i}^{\prime}f(g)] \nu(g,d\tau)
\end{eqnarray}

where $b^{i}$ and $a_{ij}$ are measurable functions on $G (1 \leq i, j \leq n)$ such that for each $g \in G$, $(a_{ij}(g))$ is a non-negative definite symmetric matrix and $\nu(g, \cdot)$ is a Borel measure on $G-\{e\}$ for which

$\int_{G-\{e\}}\left(\sum_{i=1}^{n}x_{i}(g^{-1}\tau)^{2} \wedge 1 \right)\nu(g, d\tau) < \infty$.
 It follows that ${\cal A}$ is a pseudo differential operator and by (\ref{gend}) we see  that for all $g \in G, \pi \in \G$
\begin{eqnarray} \label{Hum1}
j(g, \pi) & = & b^{i}(g)\pi(g^{-1})d\pi(X_{i})\pi(g) + a^{ij}(g)\pi(g^{-1})d\pi(X_{i})d\pi(X_{j})\pi(g) \nonumber \\
& + & \int_{G-\{e\}}[\pi(g^{-1}\tau) - I_{\pi} - x^{i}(\tau g^{-1})\pi(g^{-1})d\pi(X_{i})\pi(g)] \nu(g,d\tau)\nonumber \\
& &
\end{eqnarray}

We remark that the form (\ref{Hum}) was shown in \cite{App4} to hold for a class of Feller semigroups satisfying more general conditions (called ``Hypothesis H'' therein.)

\end{enumerate}

\section{Courr\`{e}ge-Hunt Operators and Dirichlet Forms}

Motivated by the form of the symbols appearing in the examples, we seek a general class of Markov semigroups in $L^{2}(G)$ whose infinitesimal generators are pseudo differential operators with symbol of the form
\begin{eqnarray} \label{gensy}
j(g, \pi) & = & b^{i}(g)d\pi(X_{i}) + a^{ij}(g)d\pi(X_{i})d\pi(X_{j}) \nonumber \\
& + & \int_{G-\{e\}}[\pi(g^{-1}\tau) - I_{\pi} - x^{i}(g^{-1}\tau)d\pi(X_{i})] \nu(g,d\tau),
\end{eqnarray}
for $g \in G, \pi \in \widehat{G}$. We will specify the ``characteristics'' $(b, a, \nu)$ more precisely below.

From now on we will take $G$ to be an arbitrary $n$-dimensional Lie group equipped with a {\it right-invariant} Haar measure. We define a linear operator ${\cal L}$ on $L^{2}(G)$ by the prescription\footnote{It is more convenient to work with real Hilbert spaces from now on.}
\begin{eqnarray} \label{CH}
{\cal L}f(g) & = & b^{i}(g)X_{i}f(g) + a^{ij}(g)X_{i}X_{j}f(g) \nonumber \\
& + & \int_{G}[f(\tau) - f(g) - x^{i}(g^{-1}\tau)X_{i}f(g)] \nu(g,d\tau),
\end{eqnarray}
for $f \in \Dom({\cal L}), g \in G$. We assume\footnote{We make no claim that these conditions, or those that we will impose later to study ${\cal L}^{*}$, are in any sense optimal.} that
\begin{itemize}
\item For each $1 \leq i,j,k \leq n, b^{i}$ and $a^{jk}$ are bounded measurable functions on $G$ and that for each $g \in G, (a^{jk}(g))$ is a non-negative definite symmetric matrix.

\item For each $1 \leq i \leq n, x^{i} \in C_{c}^{\infty}(G)$ and $(x_{1}, \ldots, x_{n})$ are a system of canonical co-ordinates for a given neighbourhood of the identity $U$ in $G$ that has compact closure.

\item The mapping $\nu:G \times {\cal B}(G) \rightarrow [0, \infty]$ is a {\it L\'{e}vy kernel} on $G$ so that $\nu(g, \{g\}) = 0$ for all $g \in G$, each $\nu(g, \cdot)$ defines a Borel measure on $G$, the mapping $g \rightarrow \nu(g, gE)$ is measurable for each $E \in {\cal B}(G)$ and
    $$ \sup_{g \in G}\int_{G}\left(\sum_{i=1}^{n}x_{i}(\tau)^{2} \wedge 1\right)\nu(g, gd\tau) < \infty.$$

\item There exists a L\'{e}vy measure $\rho$ on $G$ (with
$\rho(\{e\}): = 0$) such that $\nu(g, g\cdot)$ is absolutely
continuous with respect to $\rho$ for all $g \in G$.

\item $\sup_{\tau \in U^{c}}\sup_{g \in G}\lambda(g, \tau) < \infty$ where $\lambda(g, \cdot)$ is the Radon-Niko\'{y}m derivative of $\nu(g, g\cdot)$ with respect to $\rho$.

\item There exists $g^{\prime} \in G$ such that $\lambda(g^{\prime}, \tau) = \sup_{g \in G}\lambda(g, \tau)$ for all $\tau \in U$.

\end{itemize}

We call ${\cal L}$ a {\it Courr\`{e}ge-Hunt operator} as its form
is a natural extension to Lie groups of the operators on $\R^{n}$
that are considered by Courr\`{e}ge in \cite{Courr}, while when
the characteristics are constant we obtain the generators of
convolution semigroups as described by Hunt \cite{Hu}. Notice that
the form of (\ref{CH}) is also very similar to that of (\ref{Hum})
but we find it more convenient to work with left-invariant vector
fields.

Let ${\cal H}_{2}(G)$ be the Sobolev space (of order $2$) defined by
${\cal H}_{2}(G): = \{f \in L^{2}(G), X_{i}f \in L^{2}(G), X_{j}X_{k}f \in L^{2}(G), 1 \leq i,j,k \leq n\}$. It is a real Banach space under the norm
$$ |||f||||_{2} = \left(||f||_{2}^{2} + \sum_{i=1}^{n}||X_{i}f||_{2}^{2} + \sum_{j,k=1}^{n}||X_{j}X_{k}f||_{2}^{2}\right)^{\frac{1}{2}}.$$ Note that ${\cal H}_{2}(G)$ is independent of the choice of basis used to define it. We will need the fact that $C_{c}^{\infty}(G)$ is dense in ${\cal H}_{2}(G)$.

\begin{theorem} \label{Sob} ${\cal H}_{2}(G) \subseteq \Dom({\cal L})$. In fact there exists $C > 0$ such that for all $f \in {\cal H}_{2}(G)$
$$ ||{\cal L}f||_{2} \leq C||||f||||_{2}.$$ Hence ${\cal L}$ is a bounded linear operator
from ${\cal H}_{2}(G)$ to $L^{2}(G)$.
\end{theorem}

{\it Proof.}  Fix $f \in C_{c}^{\infty}(G)$ and write ${\cal L} = \cL_{1} + \cL_{2} + \cL_{3} + \cL_{4}$, where
$\cL_{1}: =  b^{i}(\cdot)X_{i}, \cL_{2}:=a^{ij}(\cdot)X_{i}X_{j}$,
$\cL_{3}f(g): = \int_{U}[f(\tau) - f(g) - x^{i}(g^{-1}\tau)X_{i}f(g)] \nu(g,d\tau)$ for each $g \in G$  and $\cL_{4}: = \cL - \cL_{1} - \cL_{2} - \cL_{3}$. We will use the elementary inequality
$$ ||\cL||^{2}_{2} \leq 4(||\cL_{1}||^{2}_{2} + ||\cL_{2}||^{2}_{2}+ ||\cL_{3}||^{2}_{2} + ||\cL_{4}||^{2}_{2}).$$
Straightforward use of the Cauchy-Schwarz inequality yields
\bean ||\cL_{1}f||_{2}^{2} & \leq & \sum_{i=1}^{n}\int_{G}b^{i}(g)^{2}\left(\sum_{j=1}^{n}X_{j}f(g)^{2}\right)dg\\
& \leq & n\max_{1 \leq i \leq n}||b_{i}(g)||^{2}_{\infty}|||f|||_{2}^{2}, \eean
\bean ||\cL_{2}f||_{2}^{2} & \leq & \sum_{i,j=1}^{n}\int_{G}a_{ij}(g)^{2}\left(\sum_{k,l=1}^{n}X_{k}X_{l}f(g)^{2}\right)dg\\
& \leq & n^{2}\max_{1 \leq i,j \leq n}||a_{i,j}(g)||^{2}_{\infty}|||f|||_{2}^{2}. \eean

We will find it convenient to make the change of variable $\tau \rightarrow g\tau$ in the integrals defining both $\cL_{3}$ and $\cL_{4}$. Then we have
\bean ||L_{4}f||^{2}_{2} & \leq & 3\int_{G}\left(\int_{U^{c}}f(g\tau)\nu(g, gd\tau)\right)^{2}dg + 3\int_{G}\left(\int_{U^{c}}f(g)\nu(g, gd\tau)\right)^{2}dg\\
& + & 3\int_{G}\left(\int_{U^{c}}\sum_{i=1}^{n}x^{i}(\tau)X_{i}f(g)\nu(g, gd\tau)\right)^{2}dg. \eean
We consider each of the three integrals separately. For the second one it is easy to
check that
$$\int_{G}\left(\int_{U^{c}}f(g)\nu(g, gd\tau)\right)^{2}dg \leq \sup_{g \in G}\nu(g, gU^{c})^{2}||f||^{2}_{2}.$$

For the first integral, we use the Cauchy-Schwarz inequality to obtain
\bean \int_{G}\left(\int_{U^{c}}f(g\tau)\nu(g, gd\tau)\right)^{2}dg & \leq & \sup_{g \in G}\nu(g, gU^{c})\int_{G}\int_{U^{c}}f(g\tau)^{2}\nu(g, gd\tau)dg\\
& = & \sup_{g \in G}\nu(g, gU^{c})\int_{G}\int_{U^{c}}f(g\tau)^{2}\lambda(g, \tau)\rho(d\tau)dg\\
& \leq & \sup_{g \in G}\nu(g, gU^{c})\sup_{\tau \in U^{c}}\sup_{g \in G}\lambda(g, \tau) \rho(U^{c})||f||^{2}_{2}, \eean
where we used Fubini's theorem and the right invariance of Haar measure to obtain the
final estimate.
Finally for the third integral we find that
\bean & & \int_{G}\left(\int_{U^{c}}\sum_{i=1}^{n}x^{i}(\tau)X_{i}f(g)\nu(g, gd\tau)\right)^{2}dg\\ & \leq & \sup_{\tau \in U^{c}}\left(\sum_{i=1}^{n}x_{i}(\tau)^{2}\right)\sup_{g \in G}\nu(g, gU^{c})^{2}\int_{G}\left(\sum_{j=1}^{n}X_{j}f(g)^{2}\right)dg\\
& \leq & n \max_{1 \leq i \leq n}||x_{i}||_{\infty}^{2}\sup_{g \in G}\nu(g, gU^{c})^{2}|||f|||^{2}_{2}. \eean
Combining the estimates for the three integrals together we can assert that there exists $K > 0$ such that
$$ ||{\cal L}_{4}f||_{2} \leq K||||f||||_{2}.$$

We now turn our attention to $\cL_{3}$ and use a Taylor expansion
as in \cite{Liao} p.13 to write for each $g \in G$,
$$ \cL_{3}f(g) = \frac{1}{2}\int_{U}\sum_{i,j=1}^{n}x_{i}(\tau)x_{j}(\tau)X_{i}X_{j}f(g\epsilon_{u}(\tau))\nu(g,gd\tau),$$
where $0 < u < 1$ and
$\epsilon_{u}(\tau):=\exp(ux^{i}(\tau)X_{i})$. By repeated
application of the Cauchy-Schwarz inequality we obtain \bean
||\cL_{3}f||^{2}_{2} & \leq &
\frac{1}{4}\int_{G}\left[\left(\int_{U}\sum_{i=1}^{n}x_{i}(\tau)^{2}\nu(g,
g d\tau)\right)\right.\times \\ & &
\left.\left(\int_{U}\sum_{i=1}^{n}x_{i}(\tau)^{2}
\sum_{j,k=1}^{n}X_{j}X_{k}f(g\epsilon_{u}(\tau))^{2}\nu(g,gd\tau)\right)\right]dg\\
& \leq & \frac{1}{4}\sup_{g \in G}\left(\int_{U}\sum_{i=1}^{n}x_{i}(\tau)^{2}\nu(g, g d\tau)\right)\times \\ & & \left( \int_{G}\int_{U}\sum_{i=1}^{n}x_{i}(\tau)^{2}\lambda(g^{\prime}, \tau)
\sum_{j,k=1}^{n}X_{j}X_{k}f(g\epsilon_{u}(\tau))^{2}\rho(d\tau)dg\right)\\
& \leq & \frac{1}{4}\sup_{g \in G}\left(\int_{U}\sum_{i=1}^{n}x_{i}(\tau)^{2}\nu(g, g d\tau)\right)\times \\ & &
\sup_{g^{\prime} \in G}\left(\int_{U}\sum_{i=1}^{n}x_{i}(\tau)^{2}\nu(g^{\prime}, g^{\prime}d\tau)\right)\int_{G}\sum_{j,k=1}^{n}X_{j}X_{k}f(g)^{2}dg\\
& \leq & \frac{1}{4}\sup_{g \in G}\left(\int_{U}\sum_{i-1}^{n}x_{i}(\tau)^{2}\nu(g,gd\tau)\right)^{2}|||f|||^{2}_{2}, \eean
where we have again utilised Fubini's theorem and the right invariance of Haar measure.


Combining together the estimates for $\cL_{1}, \cL_{2}, \cL_{3}$
and $\cL_{4}$ we obtain the existence of $C > 0$ such that $||\cL
f||_{2} \leq C|||f|||_{2}$ for all $f \in C_{c}^{\infty}(G)$ and
the required result follows by density. $\hfill \Box$

\vspace{5pt}

From now on we consider the Courr\`{e}ge-Hunt operator ${\cal L}$ as a densely defined linear operator on $L^{2}(G)$ with domain ${\cal H}_{2}(G)$.

\begin{cor}

If $G$ is compact then ${\cal L}$ is a pseudo differential operator with symbol of the form (\ref{gensy}).

\end{cor}

{\it Proof.} This is a consequence of
the fact that ${\cal M} \subseteq {\cal H}_{2}(G). \hfill \Box$

\vspace{5pt}

We would like to investigate the adjoint of ${\cal L}$ and in
order to do this we impose some additional assumptions:

\begin{itemize}

\item For all $ 1 \leq i,j,k \leq n, b_{i} \in C_{b}^{1}(G)$ and
$a_{jk} \in C_{b}^{2}(G).$

\item For all $ 1 \leq i,j,k,l,p,q \leq n, X_{l}b_{i} \in
B_{b}(G), X_{l}a_{jk} \in B_{b}(G)$ and $X_{p}X_{q}a_{jk} \in B_{b}(G)$.

\item The mapping $g \rightarrow \lambda(g, \tau) \in C^{1}(G)$
for $\rho$-almost all $\tau \in G$.

\item
$$
\int_{G}\int_{G}|x^{i}(\tau)X_{i}\lambda(g,\tau)|\rho(d\tau)dg <
\infty$$ and $c \in B_{b}(G)$ where  $c(g):=
\int_{G}x^{i}(\tau)X_{i}\lambda(g,\tau)\rho(d\tau)$ for all $g \in
G$.

\item $\lambda(g, \tau) > 0$ except on a possible set of $\gamma_{\nu}$ measure zero where
$\gamma_{\nu}(d\tau, dg):=\nu(g, gd\tau)dg$.

\item There exists $K > 0$ such that for all $g \in G, \tau \in U$
$$ |\lambda(g\tau^{-1}, \tau) - \lambda(g, \tau)| \leq K|x^{i}(\tau)\alpha_{i}(g, \tau)|,$$
where $\alpha_{i} \in B_{b}(G \times U)$ for $1 \leq i \leq n$.

\end{itemize}

We will assume that these assumptions hold for the remainder of
this paper.

We define a linear operator $R_{\rho}$ on $L^{2}(G)$ by
 the prescription
 $$ R_{\rho}f(g): = \int_{G}f(g\tau^{-1})(\lambda (g\tau^{-1}, \tau) -
 \lambda(g,\tau))\rho(d\tau).$$
for $f \in \Dom(R_{\rho}), g \in G$.
Straightforward manipulations yield

$$ R_{\rho}f(g) = \int_{G}f(\tau^{-1})\left(\frac{\lambda(\tau^{-1},\tau g)}{\lambda(g,\tau
 g)} - 1\right)\nu(g, g d\tau g).$$

  The fact that $C_{c}^{\infty}(G) \subseteq \Dom(R_{\rho})$ is demonstrated within the proof of the following theorem.

\begin{theorem} We have $C^{\infty}_{c}(G) \subseteq \Dom({\cal
L}^{*})$ and for all $f \in C^{\infty}_{c}(G), g \in G$,
\begin{eqnarray} \label{adj}
{\cal L}^{*}f(g) & = & c(g)f(g) - X^{i}(b_{i}f)(g) +
X^{i}X^{j}(a_{ij}f)(g) \nonumber \\
& = & \int_{G}\left[f(\tau^{-1})\frac{\lambda(\tau^{-1},\tau
g)}{\lambda(g,\tau g)} - f(g) + x^{i}(\tau
g)X_{i}f(g)\right]\nu(g,
g d\tau g)].\nonumber \\
& &
\end{eqnarray}
\end{theorem}

{\it Proof.} Define a linear operator $S$ on $C_{c}^{\infty}(G)$
by the action on the right hand side of (\ref{adj}). We must show
that $||Sf||_{2} < \infty$ for each $f \in C_{c}^{\infty}(G)$. We
write $S: = S^{\prime} + R_{\rho}^{(1)} + R_{\rho}^{(2)}$, where
for all $g \in G$, \bean S^{\prime}f(g) &: = & c(g)f(g) -
X^{i}(b_{i}f)(g) +
X^{i}X^{j}(a_{ij}f)(g)  \\
& = & \int_{G}\left[f(\tau^{-1}) - f(g) + x^{i}(\tau
g)X_{i}f(g)\right]\nu(g,
g d\tau g)], \eean
$R_{\rho}^{(1)}f(g): = \int_{U}f(g\tau^{-1})(\lambda (g\tau^{-1}, \tau) -
 \lambda(g,\tau))\rho(d\tau)$ and $R_{\rho}^{(2)}:=R_{\rho} - R_{\rho}^{(1)}$.
By using similar arguments to those given in the proof of Theorem \ref{Sob} we find that $||S^{\prime}f||_{2} < \infty$ and
$||R_{\rho}^{(2)}f||_{2} < \infty$. Making repeated use of the Cauchy-Schwarz inequality we find that
\bean ||R_{\rho}^{(1)}f||_{2}^{2} & = & \int_{G}\left(\int_{U}f(g \tau^{-1})(\lambda(g\tau^{-1}, \tau) - \lambda(g, \tau))\rho(d\tau)\right)^{2}dg\\
& \leq & \int_{G}\left(\int_{U}f(g \tau^{-1})^{2}|\lambda(g\tau^{-1}, \tau) - \lambda(g, \tau)|\rho(d\tau)\right)\\ & \times & \left(\int_{U}|\lambda(g\tau^{-1}, \tau) - \lambda(g, \tau)|\rho(d\tau)\right)dg\\
& \leq & K^{2}\int_{G}\left(\int_{U}f(g \tau^{-1})^{2}|x^{i}(\tau)\alpha_{i}(g, \tau)|\rho(d\tau)\right)\left(\int_{U}|x^{i}(\tau)\alpha_{i}(g, \tau)|\rho(d\tau)\right)\\
& \leq & K^{2}\sup_{g \in G, \tau \in U}\sum_{i=1}^{n}\alpha_{i}(g, \tau)^{2}\left(\int_{U}\sum_{i=1}^{n}x^{i}(\tau)^{2}\rho(d\tau)\right)^{2} ||f||^{2}_{2}\\
& < & \infty. \eean

To show that $S \subseteq {\cal L}^{*}$ its sufficient to consider the case $b_{i} = a_{ij} =
0~(1 \leq i, j \leq n)$. Let $(U_{n}, \nN)$ be a sequence in ${\cal B}(G)$ for which $U_{n} \downarrow \{e\}$.

We then find that for $f ,h \in C^{\infty}_{c}(G)$, by use of the dominated convergence theorem \bean & & \la
{\cal
L}h, f \ra \\
& = & \int_{G}\int_{G}(h(g \tau) - h(g) -
x^{i}(\tau)X_{i}h(g))f(g)\nu(g, g d\tau)dg\\
& = & \int_{G}\int_{G}(h(g \tau) - h(g) -
x^{i}(\tau)X_{i}h(g))f(g)\lambda(g,\tau)\nu(d\tau)dg\\
& = & \lim_{n \rightarrow \infty}\int_{G}\int_{G - U_{n}}(h(g
\tau) - h(g) -
x^{i}(\tau)X_{i}h(g))f(g)\lambda(g,\tau)\nu(d\tau)dg. \eean

The result follows from here by treating
each term in the integrand separately and then passing to the
limit. For example, making the change of variable $g \rightarrow g\tau^{-1}$ we have for each
$\nN$,
$$\int_{G}\int_{G - U_{n}}h(g\tau)f(g)\nu(g, g d\tau)dg  =  \int_{G}\int_{G - U_{n}}h(g)f(g\tau^{-1})\lambda(g\tau^{-1}, \tau)
\rho(d\tau)dg.$$

Using the fact that for $1 \leq i \leq n, X_{i}$ acts as a
derivation on $C_{c}^{\infty}(G)$ we obtain
$$ {\cal L}^{*}f(g)  =  c(g)f(g) + \int_{G}\left[f(g\tau^{-1})\frac{\lambda(g\tau^{-1},\tau)}{\lambda(g,\tau)} - f(g) + x^{i}(\tau)X_{i}f(g)\right]\nu(g,
g d\tau)],$$
and the result follows when we make the change of variable $\tau \rightarrow \tau g$ in the integral.$\hfill \Box $

\vspace{5pt}

We will find it convenient below to rewrite (\ref{adj}) in the form
\begin{eqnarray} \label{adj1}
{\cal L}^{*}f(g) & = & c(g)f(g) + R_{\rho}f(g) - X^{i}(b_{i}f)(g) +
X^{i}X^{j}(a_{ij}f)(g) \nonumber \\
& = & \int_{G}\left[f(\tau^{-1}) - f(g) + x^{i}(\tau
g)X_{i}f(g)\right]\nu(g,
g d\tau g)],\nonumber \\
& &
\end{eqnarray}

for each $f \in C^{\infty}_{c}(G), g \in G$. Since $C_{c}^{\infty}(G) \subseteq \Dom({\cal L})
\cap \Dom({\cal L}^{*})$, we see that both ${\cal L}$ and ${\cal
L}^{*}$ are closable in $L^{2}(G)$.

Note that ${\cal L}^{*}$ is a pseudo-differential operator in the
case when $G$ is compact. Its symbol $j^{\prime}$ is easily
calculated to be
\bean j^{\prime}(g, \pi) & = & (c(g) - (X^{i}b_{i})(g) + (X^{i}X^{j}a_{ij})(g))I_{\pi}\\
& + & (2(X^{j}a_{ij})(g) - b_{i}(g))d\pi(X^{i}) + a^{ij}(g)d\pi(X_{i})d\pi(X_{j})\\
& + & \int_{G}\left(\frac{\lambda(\tau^{-1}, \tau g)}{\lambda(g, \tau g)}\pi(g^{-1}\tau^{-1}) - I_{\pi} + x^{i}(\tau g)d\pi(X_{i})\right)\nu(g, gd\tau g), \eean
for each $g \in G, \pi \in \G$.
For the reminder of this paper we will identify ${\cal L}$ with its restriction to
$C_{c}^{\infty}(G)$.
~We can and will choose $(x_{1}, \ldots, x_{n})$ to be such that
$x_{i}(\tau^{-1}) = - x_{i}(\tau)$ for each $\tau \in G, 1 \leq i
\leq n$ (c.f. \cite{App3}, p.219.)

\begin{theorem} \label{symm} Suppose that the following conditions hold for all $g \in G$:
\begin{enumerate}
\item $b^{i}(g) = X_{j}(a^{ij})(g)$ for all $1 \leq i \leq n.$
\item 
$R_{\rho}f(g) = -c(g)f(g)$ for all $f \in C_{c}^{\infty}(G)$.
\item $\nu(g, gA) = \nu(g, gA^{-1})$ for all $A \in {\cal B}(G)$.
\end{enumerate}
Then ${\cal L}$ is symmetric and we may write
\begin{eqnarray} \label{symm1}
{\cal L}f(g)&  = &  X^{i}(a_{ij}(g)X^{j})f(g)\nonumber \\ & + &
\frac{1}{2}\int_{G}(f(g \tau) - 2f(g) + f(g\tau^{-1})\nu(g,gd\tau),
\end{eqnarray}
for all $f \in C_{c}^{\infty}(G), g \in G$.
\end{theorem}

{\it Proof.} Assume the hypotheses of the theorem hold. It is
sufficient to consider the case $a_{ij} = 0~(1 \leq i,j \leq n)$.
Making the change of variable $\tau \rightarrow g\tau^{-1}$ in
(\ref{CH}) we obtain
$$ {\cal L}f(g)
=  \int_{G}[f(g \tau^{-1}) - f(g) + x^{i}(\tau)X_{i}f(g)]
\nu(g,gd\tau).$$ We recognise the last expression as ${\cal
L}^{*}f(g)$ when we make the change of variable $\tau \rightarrow
g\tau^{-1}$ in (\ref{adj1}). Finally we obtain (\ref{symm1}) by
writing ${\cal L}(f) = \frac{1}{2}({\cal L}(f) + {\cal L}^{*}f).
\hfill \Box$

\vspace{5pt}

When ${\cal L}$ is symmetric and $G$ is compact its symbol is given by
\bean j(g, \pi) & = & (X^{i}a_{ij})(g)d\pi(X^{j}) + a_{ij}d\pi(X^{i})d\pi(X^{j})\\
& + & \frac{1}{2}\int_{G}(\pi(\tau) - 2I_{\pi} + \pi(\tau^{-1}))\nu(g, gd\tau), \eean for each $g \in G, \pi \in \G$.

Note that when $\nu \equiv 0$, then condition (1) of Theorem
\ref{symm} is both necessary and sufficient for ${\cal L}$ to be
symmetric and (\ref{symm1}) is the well-known expression for a
second order differential operator in divergence form. Suppose
that condition (1) of Theorem \ref{symm} holds and that for each
$g \in G, A \in {\cal B}(G), \nu(g,A) = \rho(g^{-1}A)$ where
$\rho$ is (as above) a L\'{e}vy measure on $G$. In this case
condition (3) of Theorem \ref{symm} is just the requirement that
$\rho$ is a symmetric measure and condition 2 holds automatically.
We then have ${\cal L} = {\cal L}_{D} + {\cal L}_{J}$ where ${\cal
L}_{D}$ is the infinitesimal generator of a diffusion process on
$G$ (in divergence form) and ${\cal L}_{J}$ generates a pure jump
L\'{e}vy process on $G$. It may be that this is the only
possibility (however condition (2) in Theorem \ref{symm} is very strong) - indeed for future work on symmetric Courr\'{e}ge-Hunt operators it may be more fruitful to work in $L^{2}(G, \mu)$ where $\mu$ is an infinitesimal invariant
measure for the linear operator ${\cal L}$ (see \cite{ARW}), but in that case the representation using pseudo differential operators will require further development.

\vspace{5pt}

We assume that ${\cal L}$ is symmetric and define the symmetric
bilinear form
$$ {\cal E}(f_{1}, f_{2}) = -\la {\cal L}f_{1}, f_{2} \ra,$$
for $f_{1}, f_{2} \in C_{c}^{\infty}(G)$.

\begin{cor} The form ${\cal E}$ has an extension which is a Dirichlet form on $L^{2}(G)$. It has the Beurling-Deny representation
\begin{eqnarray} \label{BD} {\cal E}(f_{1},f_{2}) & = &
\int_{G}a^{ij}(g)(X_{i}f_{1})(g)(X_{j}f_{2})(g)dg \nonumber \\ & + & \frac{1}{2}\int_{(G
\times G) - D}(f_{1}(\tau) - f_{1}(g))(f_{2}(\tau) - f_{2}(g)J(d\tau,
dg), \nonumber\\ & & \end{eqnarray} for $f_{1}, f_{2} \in C_{c}^{\infty}(G)$ where $D: = \{(g, g), g \in G\}$ and
$J(d\tau, dg): = \nu(g, d\tau)dg$.
\end{cor}

{\it Proof.} The representation (\ref{BD}) is obtained by standard manipulations (see e.g. Theorem 2.4 in \cite{App3} and Proposition 2.1 in \cite{Kun} for the case where ${\cal L}$ is the generator of a convolution semigroup.) It then follows that ${\cal E}$ is Markovian by the argument given on page 7 of \cite{FOT}. The form is closeable by standard arguments. Hence by Theorem 3.1.1 of \cite{FOT}, p.98 its smallest closed extension is a Dirichlet form. $\hfill \Box$

\vspace{5pt}

Some interesting relations between closability of ${\cal E}$ and the behaviour of certain negligible sets can be deduced from results in \cite{AS}.

\vspace{5pt}

{\it Acknowledgements.} I would like to thank Michael Ruzhansky
for several helpful comments on an early draft of this paper.


\begin{thebibliography}{99}

\bibitem{ARW} S.Albeverio, B.R\"{u}diger, J-L.Wu, Invariant
measures and symmetry property of L\'{e}vy type operators, {\it
Potential Anal.} {\bf 13}, 147-68 (2000)

\bibitem{AS} S.Albeverio, S.Song, Closability and resolvent of
Dirichlet forms perturbed by jumps, {\it Potential Anal.} {\bf
2}, 115-30 (1993)

\bibitem{App1} D.Applebaum, Operator-valued stochastic differential equations
arising from unitary group representations, {\it Journal of
Theoretical Probability} {\bf 14}, 61-76 (2001)

\bibitem{App4} D.Applebaum, On a class of non-translation invariant Feller semigroups on
Lie groups, {\it Potential Anal.} {\bf 16}, 103-14 (2002)

\bibitem{App2} D.Applebaum,  Infinitely divisible central probability measures on compact
Lie groups - regularity, semigroups and transition kernels, to
appear in {\it Annals of Prob.} (2010)

\bibitem{Abook} D.Applebaum {\it L\'{e}vy Processes and Stochastic
Calculus} (second edition), Cambridge University Press (2009)

\bibitem{App3} D.Applebaum, Some $L^{2}$ properties of semigroups
of measures on Lie groups, {\it Semigroup Forum} {\bf 79}, 217-28
(2009)

\bibitem{ATa2} D.Applebaum, F.Tang, Stochastic flows of
diffeomorphisms on manifolds driven by infinite-dimensional
semimartingales with jumps, {\it Stoch. Proc. Appl.} {\bf 92},
219--36 (2002).

\bibitem{BHT} F.Baudoin, M.Hairer, J.Teichmann, Ornstein-Uhlenbeck processes on Lie groups,
 {\it J. Funct. Anal.} {\bf 255}, 877-90 (2008)

\bibitem{BCP} J.M.Bony, P.Courr\`{e}ge, P.Prioret, Semi-groupes de Feller sur
une vari\'{e}t\'{e} a bord compacte et probl\`{e}mes aux limites
int\'{e}gro-diff\'{e}rentiels du second-ordre donnant lieu au
principe du maximum, {\it Ann. Inst. Fourier, Grenoble} {\bf 18},
369-521 (1968)

\bibitem{Courr} P.Courr\`{e}ge, Sur la forme int\'{e}gro-diff\'{e}rentielle des
op\'{e}rateurs de $C_{k}^{\infty}$ dans $C$ satifaisant au
principe du maximum, {\it S\'{e}m. Th\'{e}orie du Potential}
expos\'{e} {\bf 2},  (1965/66) 38pp.

\bibitem{Courr1} P.Courr\`{e}ge, Sur la forme int\'{e}gro-diff\'{e}rentielle du g\'{e}n\'{e}rateur infinit\'{e}simal d'un semi-groupe de Feller
sur une vari\'{e}t\'{e}, {\it S\'{e}m. Th\'{e}orie du Potential}
expos\'{e} {\bf 3},  (1965/66) 48pp.


\bibitem{EK} S.N.Ethier, T.G.Kurtz, {\it Markov Processes,
Characterisation and Convergence}, Wiley (1986).

\bibitem{Far} J.Faraut, {\it Analysis on Lie Groups}, Cambridge
University Press (2008)

\bibitem{FOT} M. Fukushima, Y. Oshima, M. Takeda, {\it Dirichlet
Forms and Symmetric Markov Processes},   de Gruyter (1994).

\bibitem{Hey} H.Heyer, L'analyse de Fourier
non-commutative et applications \`{a} la th\'{e}orie des
probabilit\'{e}s, {\it Ann. Inst. Henri Poincar\'{e} (Prob.
Stat.)} {\bf 4}(1968), 143-68

\bibitem{Hey1} H.Heyer, Infinitely divisible probability measures
on compact groups, in Lectures on Operator Algebras pp.55-249,
Lecture Notes in Math. vol {\bf 247} Springer Berlin, Heidelberg,
New York (1972)

\bibitem{He1} H.Heyer, {\it Probability Measures on Locally Compact Groups}, Springer-Verlag,
Berlin-Heidelberg (1977)

\bibitem{Hoh1} W.Hoh, The martingale problem for a class of pseudo differential operators, {\it Math. Ann.} {\bf 300}, 121-47 (1994)

\bibitem{Hoh2} W.Hoh, Pseudo differential operators generating Markov processes, Habilitationsschrift, Bielefeld (1998)

\bibitem{HJ} W.Hoh, N.Jacob, Some Dirichlet forms generated by pseudo differential operators, {\it Bull. Sc. Math.} 2 s\'{e}rie, {\bf 116}, 383-98 (1992)

\bibitem{Hor} L.H\"{o}rmander, {\it The Analysis of Linear Partial
Differential Operators III},{ \it Grundlehren der Mathematischen
Wissenshaften} {\bf 274}, Springer-Verlag Berlin Heidelberg (1985)

\bibitem{Hu} G.A.Hunt, Semigroups of measures on Lie groups, {\it Trans. Amer.
Math. Soc.} {\bf 81}, 264-93 (1956)

\bibitem{Jacold} N.Jacob, Translation invariant pseudo
differential operators on compact abelian groups, unpublished
(1987)

\bibitem{Jac0} N.Jacob, A class of Feller semigroups generated by pseudo differential
operators, {\it Math. Z.} {\bf 215}, 151-66 (1994)

\bibitem{Jac1} N.Jacob, {\it Pseudo-Differential Operators and Markov Processes}, Akademie-Verlag, Mathematical Research 94 (1996).

\bibitem{Jac2} N.Jacob, {\it Pseudo Differential Operators and Markov Processes:
 1, Fourier Analysis and Semigroups}, World Scientific (2001).

\bibitem{Jac3} N.Jacob, {\it Pseudo Differential Operators and Markov Processes:
 2, Generators and Their Potential Theory}, World Scientific (2002).

\bibitem{Jac4} N.Jacob, {\it Pseudo Differential Operators and Markov Processes:
 3, Markov Processes and Applications}, World Scientific (2005).

 \bibitem{Kato} T.Kato, {\it Perturbation Theory for Linear
Operators} (second edition), Springer-Verlag (1995).

\bibitem{Kom} T. Komatsu, Markov processes associated with certain
integro-differential operators, {\it Osaka J. Math.} {\bf 10},
271-305 (1973).

\bibitem{Kun} H.Kunita, Analyticity and injectivity of convolution
semigroups on Lie groups, {\it J.Funct. Anal.} {\bf 165}, 80-100
(1999)

\bibitem{Liao} M.Liao, {\it L\'{e}vy Processes in Lie
Groups}, Cambridge University Press (2004)

\bibitem{LoNg} J.T-H.Lo, S-K.Ng, Characterizing Fourier series
representations of probability distributions on compact Lie
groups, {\it Siam J. Appl. Math.} {\bf 48} 222-8 (1988)

\bibitem{Phi} R.S.Phillips, On the generation of semigroups of linear operators, {\it Pacific J.Math.}
{\bf 2}, 343-69 (1952)

\bibitem{RT} M.Ruzhansky, V.Turunen, {\it Pseudo-differential
Operators and Symmetries: Background Analysis and Advanced
Topics}, Birkh\"{a}user, Basel (2010)

\bibitem{RTW} M.Ruzhansky, V.Turunen, J.Wirth, H\"{o}rmander class of pseudo-differential operators on compact Lie groups and global hypoellipticity, preprint   arXiv:1004.4396v1 (2010)

\bibitem{Sa} K.-I.Sato, {\it L\'{e}vy Processes and Infinite
Divisibility}, Cambridge University Press (1999)

\bibitem{Sch} R.L.Schilling, Conservativeness of semigroups
generated by pseudo differential operators, {\it Potential
Analysis} {\bf 9}, 91-104 (1998)

\bibitem{ScSc} R.L.Schilling, A.Schnurr, The symbol associated with the solution of a stochastic differential equation, {\it Electronic J. Prob.} {\bf 15}, 1369-1393 (2010)

\bibitem{Sieb} E.Siebert, Fourier analysis and limit
theorems for convolution semigroups on a locally compact group,
{\it Advances in Math.} {\bf 39}, 111-54 (1981)

\bibitem{Stro} D. Stroock, Diffusion processes associated with
L\'{e}vy generators, {\it Z. Wahrsch. verw. Geb.} {\bf 32}, 209-44 (1975).


\bibitem{Tay} M.Taylor, {\it Pseudo Differential Operators},
Lecture Notes in Mathematics {\bf 416}, Springer-Verlag Berlin
Heidelberg (1974)







\end{thebibliography}
\end{document}